\documentclass[a5paper]{book}
\usepackage[english]{babel}
\usepackage{hyperref,amsthm,amsmath,amssymb,color,multirow,graphicx}
\usepackage{inputenc}
\usepackage{wrapfig, float,array}

\newtheorem{thm}{Theorem}[section]

\theoremstyle{definition}
\newtheorem{defn}[thm]{Definition}

\numberwithin{equation}{section}

\newtheorem*{thm*}{Теорема}
\newtheorem*{lem*}{Лемма}

\newcommand{\be}{\begin{equation}}

\newcommand{\ee}{\end{equation}}

\begin{document}

\sloppy


\sloppy

\begin{center}
\textbf{\large Automorphisms of some solvable Leibniz algebras}\\

\textbf{I.A.Karimjanov, S.M.Umrzaqov}
\end{center}

{\small \textbf{Abstract.} In this paper, we describe the automorphisms of solvable Leibniz algebras with null-filiform nilradical. Moreover we describe the automorphisms of solvable Leibniz algebras with naturally graded non-Lie filiform nilradicals, whose the dimension of complementary space is maximal.
\\


\\

\makeatletter
\renewcommand{\@evenhead}{\vbox{\thepage \hfil {\it I.A.Karimjanov, S.M.Umrzaqov}   \hrule }}
\renewcommand{\@oddhead}{\vbox{\hfill
{\it Isomorphisms of some solvable Leibniz algebras}\hfill
\thepage \hrule}} \makeatother

\label{firstpage}

\section{Introduction}

\qquad Leibniz algebras were defined by J.-L.Lodey \cite{loday} at the beginning of the 90s of
the past century. Leibniz algebras generalize Lie algebras in natural way. In fact, many results of the theory of Lie algebras are extended to Leibniz algebras case. During the last 30 years the theory of Leibniz algebras has been actively investigated and numerous papers have been devoted to the study of these algebras. The classification problem of finite-dimensional Leibniz algebras is very difficult problem. But a lot of main structural results were obtained. For instance, the classic Levi-Malcev decomposition of Leibniz algebras were obtained by D.W. Barnes \cite{barn}, the classical results on Cartan subalgebras and Engel's theorem are established in Leibniz algebras case, proved Mubarakzjanov theorem  for the Leibniz algebras case \cite{cas1} which says that the dimension of the splittable algebra does not exceed the number of nil-independent derivations of
the nilradical.

\quad Automorphism is a fundamental notion in mathematics. Automorphisms play a prominent role in algebra. In the context of abstract algebra, a mathematical object is an algebraic structure such as a group, ring, or vector space. An automorphism is simply a bijective homomorphism of an object with itself. There are a huge applications of automorphisms for theory of group, number theory and etc. There are many generalizations of automorphism. The important generalization of automorphism are local and 2-local automorphisms. The last years many authors published several works about local and 2-local automorphisms of some varieties.

\quad We obtain some very basic results concerning automorphism forms for solvable Leibniz algbras with null-filiform nilradical and naturally graded non-Lie filiform nilradical of maximal dimension.

\section{Preliminaries}

All algebras considered are supposed to be over the field of complex numbers $\mathbb{C}$.

\begin{defn} An algebra $L$ is called a Leibniz algebra if for any $x,y,z\in L$, the Leibniz identity
		\[[[x,y],z]=[[x,z],y]+[x,[y,z]]\] is satisfied, where $[-,-]$ is the multiplication in $L$.
	\end{defn}

For a Leibniz algebra $L$ we consider the following derived and lower central series:
	\begin{align*}
		& \text{(i)}   &L^{[1]}= & \ L, \  & L^{[n+1]}= & \ [L^{[n]},L^{[n]}],  & n>1; \\
		& \text{(ii)}  & L^1= &  \ L, \ & L^{n+1}=&  \ [L^n,L],  & n>1.
	\end{align*}

	\begin{defn} An algebra $L$ is called solvable (nilpotent) if there exists $s\in \mathbb{N}$ \ ($k\in \mathbb{N}$, respectively) such that $L^{[s]}=0$ \ ($L^k=0$, respectively).
The minimal number $s$ (respectively, $k$) with such property is called index of
solvability (respectively, of nilpotency) of the algebra $L.$
	\end{defn}

Evidently, the index of nilpotency of  an $n$-dimensional algebra is not greater than $n+1.$

\begin{defn} An $n$-dimensional Leibniz algebra is called null-filiform if $\dim L^i=n+1-i, \ 1\leq i \leq n+1.$
\end{defn}

Actually, a nilpotent Leibniz algebra is null-filiform if it is a one-generated algebra. Note, that this notion has no sense in Lie algebras case, because they are at least two-generated.

\begin{defn} A Leibniz algebra $L$ is said to be filiform if
$\dim L^i=n-i$, for $2\leq i \leq n$, where $n=\dim L$.
\end{defn}

\begin{defn} Given a filiform Leibniz algebra $L,$ put
$L_i=L^i/L^{i+1}, \ 1 \leq i\leq n-1,$ and $gr L = L_1 \oplus
L_2\oplus\dots L_{n-1}.$ Then $[L_i,L_j]\subseteq L_{i+j}$ and we
obtain the graded algebra $gr L$. If $gr L$ and $L$ are isomorphic,
denoted by $gr L\cong L,$ we say that the algebra $L$ is naturally
graded.
\end{defn}

\begin{defn} The  (unique) maximal nilpotent ideal of a Leibniz algebra is called the nilradical of the algebra.
\end{defn}

Let $R$ be a solvable Leibniz algebra. Then it can be decomposed into the form $R=N \oplus Q$, where $N$ is the  nilradical and $Q$ is the complementary vector space. Since the square of a solvable algebra is a nilpotent ideal and the finite sum of nilpotent ideals is a nilpotent ideal too, then the ideal $R^2$ is nilpotent, i.e. $R^2\subseteq N$ and consequently, $Q^2\subseteq N$.

Now, we present the classification results for arbitrary dimensional solvable Leibniz algebras with null-filiform and naturally graded non-Lie filiform nilradicals, whose the dimension of complementary space is maximal.

\begin{thm} \cite{cas1} Let $R_0$ be a solvable Leibniz algebra with null-filiform nilradical. Then there exists a basis $\{e_0, e_1, e_2, \dots, e_n\}$ of the algebra $R_0$ such that the multiplication table of $R_0$ with respect to this basis has the following form:
\[R_0: \left\{ \begin{aligned}
{}[e_i,e_1] & =e_{i+1}, && 0\leq i\leq n-1,\\
[e_i,e_0] & =-ie_i, && 1\leq i\leq n.
\end{aligned}\right.\]
\end{thm}

\begin{thm}\cite{cas2,kam} An arbitrary $(n+2)$-dimensional solvable Leibniz algebra with n-dimensional naturally graded non-Lie filiform nilradical is isomorphic to one of the following non-isomorphic algebras:
\[R_1: \left\{\begin{array}{ll}
[e_i,e_1]=e_{i+1}, &  2\leq i \leq {n-1},  \\[1mm]
[e_1,x]=-[x,e_1]=e_1,      &               \\[1mm]
[e_i,x]=(i-1)e_i, &  2\leq i \leq n, \\[1mm]
[e_i,y]=e_i,       &  2\leq i\leq n,     \\[1mm]
\end{array}\right.\]
\[ R_2: \left\{\begin{array}{ll}
[e_1,e_1]=e_3, & \\[1mm]
[e_i,e_1]=e_{i+1}, & 3\leq i\leq n-1,\\[1mm]
[e_1,x]=-[x,e_1]=e_1, & \\[1mm]
[e_i,x]=(i-1)e_i, & 3\leq i\leq n, \\[1mm]
[e_2,y]=-[y,e_2]=e_2, &
\end{array}\right.\]
\[ R_3: \left\{\begin{array}{ll}
[e_1,e_1]=e_3, & \\[1mm]
[e_i,e_1]=e_{i+1}, & 3\leq i\leq n-1,\\[1mm]
[e_1,x]=-[x,e_1]=e_1, & \\[1mm]
[e_i,x]=(i-1)e_i, & 3\leq i\leq n, \\[1mm]
[e_2,y]=e_2 &
\end{array}\right.\]
where $\{e_1, \dots, e_n, x, y\}$ is a basis of the algebra.
\end{thm}

\section{Main Result}

Now we shall give the main result concerning automorphisms of solvable Leibniz algebra with null-filiform nilradicals.

\begin{thm} \label{1} A linear map $\varphi:R_0\to R_0$ is a automorphism if and only if when $\varphi$ has the
following form:
\[\varphi(e_i)=\sum\limits_{j=i}^n \frac{\alpha^{j-i} \beta^{i}}{(j-i)!} e_j, \quad   0\leq i \leq n, \]
where \(\beta \neq  0\).
\end{thm}
\begin{proof} We denote an automorphism $\varphi$ acts as follows:
\[\varphi(e_i)=\sum\limits_{j=0}^na_{j,i}e_j, \quad 0\leq i\leq1, \qquad \varphi(e_k)=[\varphi(e_{k-1}),\varphi(e_1)],\quad 2\leq k \leq n.\]

Considering the equalities
\[\begin{array}{c}
0=\varphi([e_0,e_0])=[\varphi(e_0),\varphi(e_0)]=[\sum\limits_{i=0}^na_{i,0}e_i,\sum\limits_{i=0}^na_{i,0}e_i]=\\[1mm]
\sum\limits_{i=0}^na_{0,0}a_{i,0}[e_i,e_0]+\sum\limits_{i=0}^na_{1,0}a_{i,0}[e_i,e_1]=\sum\limits_{i=2}^n(a_{1,0}a_{i-1,0}-ia_{0,0}a_{i,0})e_i\\[1mm]
\end{array}\]
we obtain the next recurrent expression:
\[a_{0,0}a_{i,0}=\frac{a_{1,0}a_{i-1,0}}{i}, \quad 2\leq i\leq n\]
which immediately implies
\begin{equation}\label{null1}
a_{0,0}^{i-1}a_{i,0}=\frac{a_{1,0}^{i}}{i!}, \quad 2\leq i\leq n.
\end{equation}

Consider
\[\begin{array}{c}
\sum\limits_{i=0}^na_{i,1}e_i=\varphi(e_1)=\varphi([e_0,e_1])=[\varphi(e_0),\varphi(e_1)]=[\sum\limits_{i=0}^na_{i,0}e_i,\sum\limits_{i=0}^na_{i,1}e_i]=\\
\sum\limits_{i=0}^na_{0,1}a_{i,0}[e_i,e_0]+\sum\limits_{i=0}^na_{1,1}a_{i,0}[e_i,e_1]=\sum\limits_{i=1}^n(a_{1,1}a_{i-1,0}-ia_{0,1}a_{i,0})e_i.
\end{array}\]
Comparing the coefficients at the basis elements we derive: $a_{0,1}=0, a_{0,0}=1$ and
\[a_{i,1}=a_{i-1,0}a_{1,1}, \quad 2\leq i\leq n.\]
By using the equality (\ref{null1}) we obtain
\[a_{i,1}=\frac{a_{1,0}^{i-1}a_{1,1}}{(i-1)!}, \quad 2\leq i\leq n.\]

So, we have
\[\varphi(e_0)=\sum\limits_{i=0}^n\frac{a_{1,0}^i}{i!}e_i, \quad \varphi(e_1)=\sum\limits_{i=1}^n\frac{a_{1,0}^{i-1}a_{1,1}}{(i-1)!}e_i.\]

Now we shall prove the following equalities by an induction on $i$:
\begin{equation}\label{null2}
\varphi(e_i)=\sum\limits_{j=i}^n \frac{a_{1,0}^{j-i} a_{1,1}^{i}}{(j-i)!} e_j, \quad   0\leq i \leq n.
 \end{equation}

 Obviously, the equality holds for $i = 0,1$. Let us assume that the equality holds
for $2 < i < n$, and we shall prove it for $i + 1$:
\[\begin{array}{c}
\varphi(e_{i+1})=[\varphi(e_i),\varphi(e_1)]=[\sum\limits_{j=i}^n \frac{a_{1,0}^{j-i} a_{1,1}^i}{(j-i)!} e_j,\sum\limits_{k=1}^n\frac{a_{1,0}^{k-1}a_{1,1}}{(k-1)!}e_k]=\\[1mm]
\sum\limits_{j=i}^n \frac{a_{1,0}^{j-i} a_{1,1}^{i+1}}{(j-i)!} [e_j,e_1]=\sum\limits_{j=i+1}^n\frac{a_{1,0}^{j-(i+1)}a_{1,1}^{i+1}}{(j-(i+1))!}e_j.
\end{array}\]
so, the induction proves the equalities (\ref{null2}) for any $i, 0 \leq i \leq n$.

By denoting $(a_{1,0},a_{1,1})=(\alpha,\beta)$ we complete necessity condition of the proof.

Now, we shall prove the sufficiency condition of the theorem.

Let us $u,v\in R_0$ and denote by
\[u=\sum\limits_{i=0}^n\lambda_ie_i,
\quad v=\sum\limits_{i=0}^n\mu_ie_i.\]

Consider
\[\begin{array}{c}
\varphi([u,v])=\varphi\big(\big[\sum\limits_{i=0}^n\lambda_ie_i,\sum\limits_{i=0}^n\mu_ie_i \big]\big)=\\
\varphi( \mu_1 \sum\limits_{i=1}^{n}\lambda_{i-1}e_i -\mu_0 \sum\limits_{i=1}^n i \lambda_i e_i)=\varphi(\sum\limits_{i=1}^n(\mu_1 \lambda_{i-1} -i\mu_0 \lambda_i)e_i)=\\
=\sum\limits_{i=1}^n(\mu_1 \lambda_{i-1} -i\mu_0 \lambda_i)\sum\limits_{j=i}^n \frac{\alpha^{j-i}\beta^i}{(j-1)!}e_j=\\
=\sum\limits_{i=1}^{n}\sum\limits_{j=1}^{i}(\mu_1 \lambda_{j-1}-j \mu_0 \lambda_j)\frac{\alpha^{i-j} \beta^j}{(i-j)!}e_i.
\end{array}\]

On the other hand,
\[\begin{array}{c}
[\varphi(u),\varphi(v)]=[\varphi(\sum\limits_{i=0}^n\lambda_ie_i),\varphi(\sum\limits_{i=0}^n\mu_ie_i)]=\\[1mm]
=[\sum\limits_{i=0}^n\lambda_i \sum\limits_{j=i}^n \frac{\alpha^{j-i} \beta^i}{(j-i)!}e_j,\sum\limits_{i=0}^n\mu_i \sum\limits_{j=i}^n \frac{\alpha^{j-i} \beta^i}{(j-i)!}e_j]=\\[1mm]
[\sum\limits_{i=0}^n \sum\limits_{j=0}^i \frac{\lambda_j \alpha^{i-j} \beta^j}{(i-j)!} e_i,\sum\limits_{i=0}^n \sum\limits_{j=0}^i \frac{\mu_j \alpha^{i-j} \beta^j}{(i-j)!} e_i ] =\\[1mm]
=(\sum\limits_{i=0}^{n-1} \sum\limits_{j=0}^i \frac{\lambda_j \alpha^{i-j} \beta^j}{(i-j)!}) (\sum\limits_{j=0}^1 \frac{\mu_j \alpha^{1-j}\beta^j}{(1-j)!}) [e_i,e_1]+(\sum\limits_{i=1}^{n} \sum\limits_{j=0}^i \frac{\lambda_j \alpha^{i-j} \beta^j}{(i-j)!}) (\mu_0)[e_i, e_0]=\\[1mm]
=(\mu_0 \alpha+\mu_1 \beta)\sum\limits_{i=0}^{n-1} \sum\limits_{j=0}^i \frac{\lambda_j \alpha^{i-j} \beta^j}{(i-j)!} e_{i+1}-\mu_0 \sum\limits_{i=1}^n \sum\limits_{j=0}^i \frac{i \lambda_j \alpha^{i-j} \beta^j}{(i-j)!} e_i =\\[1mm]
(\mu_0 \alpha+\mu_1 \beta)\sum\limits_{i=1}^{n} \sum\limits_{j=0}^{i-1} \frac{\lambda_j \alpha^{i-j-1} \beta^j}{(i-j-1)!} e_i-\mu_0 \sum\limits_{i=1}^n \sum\limits_{j=0}^i \frac{i \lambda_j \alpha^{i-j} \beta^j}{(i-j)!} e_i =\\[1mm]
=\sum\limits_{i=1}^{n}((\mu_0 \alpha+\mu_1 \beta)\sum\limits_{j=0}^{i-1} \frac{\lambda_j \alpha^{i-j-1} \beta^j}{(i-j-1)!}-\mu_0 \sum\limits_{j=0}^i \frac{i \lambda_j \alpha^{i-j} \beta^j}{(i-j)!})e_i=\\
=\sum\limits_{i=1}^{n}(\mu_0 \sum\limits_{j=1}^{i}\frac{-j \lambda_j \alpha^{i-j} \beta^j}{(i-j)!}+\mu_1 \sum\limits_{j=1}^{i}\frac{\lambda_{j-1} \alpha^{i-j} \beta^j}{(i-j)!})e_i= \\
=\sum\limits_{i=1}^{n}\sum\limits_{j=1}^{i}(\mu_1 \lambda_{j-1}-j \mu_0 \lambda_j)\frac{\alpha^{i-j} \beta^j}{(i-j)!}e_i. \\
\end{array}\]

Comparing coefficients at the basis elements we obtain that
\[\varphi([u,v])=[\varphi(u),\varphi(v)]\]
and we complete the proof of theorem.
\end{proof}

Now we present the main result concerning automorphisms on solvable Leibniz algebras
with naturally graded non-Lie filiform nilradicals, whose the dimension of complementary space is maximal.

\begin{thm} A linear map $\varphi:R_1\to R_1$ is a automorphism if and only if when $\varphi$ has the
following form:
\[\begin{array}{ll}
     \varphi(e_1)=\alpha e_1, & \\
      \varphi(e_i)=\sum\limits_{j=i}^n \frac{(-1)^{j-i}\alpha^{i-2} \beta\gamma^{j-i}}{(j-i)!} e_j,&  2\leq i\leq n,\\
       \varphi(x)=\gamma e_1+x, & \\
        \varphi(y)=y, & \\
     \end{array}\]
 where \(\alpha \beta \neq0\).
\end{thm}

\begin{proof} We set
\[\begin{array}{ll}
    \varphi(e_i)=\sum\limits_{j=1}^na_{j,i}e_j+a_{n+1,i}x+a_{n+2,i}y, & 1\leq i\leq 2,\\
      \varphi(e_k)=\varphi([e_{k-1},e_1])=[\varphi(e_{k-1}),\varphi(e_1)],&  3\leq k\leq n,\\
      \varphi(x)=\sum\limits_{j=1}^na_{j,n+1}e_j+a_{n+1,n+1}x+a_{n+2,n+1}y,  & \\
        \varphi(y)=\sum\limits_{j=1}^na_{j,n+2}e_j+a_{n+1,n+2}x+a_{n+2,n+2}y. & \\
     \end{array}\]

From the equalities
\[\varphi(e_1)=\varphi([e_1,x])=[\varphi(e_1),\varphi(x)], \quad  \varphi(e_i)=\varphi([e_i,y])=[\varphi(e_i),\varphi(y)], \quad 2\leq i\leq n,\]
we derive that $\varphi(e_i)\in span(e_1,\dots,e_n)$ for $1\leq i\leq n$.

From the multiplication table it is easy to see that $e_i\in Ann_r(R_1)$ for $2\leq i\leq n$. So, without loss of generality, we can assume that $a_{1,1}\neq0$, otherwise the matrix of automorphism is singular.

Considering the equalities:
\[0=\varphi([e_1,e_1])=[\varphi(e_1),\varphi(e_1)]=[\sum\limits_{i=1}^na_{i,1}e_i,\sum\limits_{i=1}^na_{i,1}e_i]=a_{1,1}\sum\limits_{i=3}^na_{i-1,1}e_i,\]
we can deduce that $a_{i,1}=0$ for $2\leq i\leq n-1$ or
\[\varphi(e_1)=a_{1,1}e_1+a_{n,1}e_n.\]

From the following chain of equalities:
\[\begin{array}{c}
0=\varphi([y,e_1])=[\varphi(y),\varphi(e_1)]=\\[1mm]
[\sum\limits_{j=1}^na_{j,n+2}e_j+a_{n+1,n+2}x+a_{n+2,n+2}y,a_{1,1}e_1+a_{n,1}e_n]=\\[1mm]
a_{1,1}(-a_{n+1,n+2}e_1+\sum\limits_{i=3}^na_{i-1,n+2}e_i),
\end{array}\]
we obtain $a_{n+1,n+2}=a_{i,n+2}=0, \ 2\leq i\leq n-1$ or
\[\varphi(y)=a_{1,n+2}e_1+a_{n,n+2}e_n+a_{n+2,n+2}y.\]

Consider
\[\begin{array}{c}
a_{1,1}e_1+a_{n,1}e_n=\varphi(e_1)=-[\varphi([x,e_1])]=-[\varphi(x),\varphi(e_1)]=\\[1mm]
[\sum\limits_{j=1}^na_{j,n+1}e_j+a_{n+1,n+1}x+a_{n+2,n+1}y,a_{1,1}e_1+a_{n,1}e_n]=\\[1mm]
a_{1,1}(a_{n+1,n+1}e_1-\sum\limits_{i=3}^na_{i-1,n+1}e_i).
\end{array}\]

Comparing coefficients at the basis elements in right and left sides, we
obtain the following relations:
\[a_{n+1,n+1}=1, \quad a_{n,1}=-a_{n-1,n+1}a_{1,1}, \quad a_{i,n+1}=0, \quad 2\leq i\leq n-2.\]

Following we have
\[\varphi(x)=a_{1,n+1}e_1+a_{n-1,n+1}e_{n-1}+a_{n,n+1}e_n+x+a_{n+2,n+1}y.\]

The equalities
\[\begin{array}{c}
0=\varphi([x,e_2])=[\varphi(x),\varphi(e_2)]=\\
=[a_{1,n+1}e_1+a_{n-1,n+1}e_{n-1}+a_{n,n+1}e_n+x+a_{n+2,n+1}y,\sum\limits_{i=1}^na_{i,2}e_i]=\\
-a_{1,2}e_1+a_{1,2}a_{n-1,n+1}e_n,
\end{array}\]
yield $a_{1,2}=0$, and it follows:
\[\varphi(e_2)=\sum\limits_{i=2}^na_{i,2}e_i.\]

From the next equalities, we have
\[
\begin{array}{lllll}
\text{Equality }& & & & \text{ Constraint }\\[1mm]
\hline \hline\\[1mm]
\varphi([y,x])&=&[\varphi(y),\varphi(x)]&\Longrightarrow &a_{1,n+2}=0, \\[1mm]
\varphi([e_2,y])&=&[\varphi(e_2),\varphi(y)]&\Longrightarrow &a_{n+2,n+2}=1, \\[1mm]
\varphi([y,y])&=&[\varphi(y),\varphi(y)]&\Longrightarrow &a_{n,n+2}=0, \\[1mm]
\varphi([e_1,y])&=&[\varphi(e_1),\varphi(y)]&\Longrightarrow &a_{n,1}=0, \\[1mm]
\varphi([x,y])&=&[\varphi(x),\varphi(y)]&\Longrightarrow &a_{n-1,n+1}=a_{n,n+1}=0.
\end{array}
\]

Now we consider the following equlities
\[\begin{array}{c} \sum\limits_{i=2}^na_{i,2}e_i=\varphi(e_2)=\varphi([e_2,x])=[\varphi(e_2),\varphi(x)]=\\[1mm]
[\sum\limits_{i=2}^na_{i,2}e_i,a_{1,n+1}e_1+x+a_{n+2,n+1}y]=\\\\[1mm]
=\sum\limits_{i=3}^na_{i-1,2}a_{1,n+1}e_i+\sum\limits_{i=2}^n(i-1)a_{i,2}e_i+\sum\limits_{i=2}^na_{i,2}a_{n+2,n+1}e_i=\\\\[1mm]
=a_{2,2}(1+a_{n+2,n+1})e_2+\sum\limits_{i=3}^n(a_{i-1,2}a_{1,n+1}+a_{i,2}(i-1+a_{n+2,n+1}))e_i
\end{array}\]
By comparison of coefficients at the elements of the basis we obtain that $a_{n+2,n+1}=0$ and the next recurrent expression:
\[a_{i,2}=-\frac{a_{1,n+1}a_{i-1,2}}{i-2}, \quad 3\leq i\leq n.\]
From the previous recurrent expressions we can deduce that
\[a_{i,2}=(-1)^i\frac{a_{2,2}a_{1,n+1}^{i-2}}{(i-2)!}, \quad 3\leq i\leq n.\]
Consequently,
\[\varphi(e_2)=\sum\limits_{i=2}^n(-1)^i\frac{a_{2,2}a_{1,n+1}^{i-2}}{(i-2)!}e_i.\]

Now we shall prove the following equalities by an induction on $i$:
\begin{equation}\label{r1}
 \varphi(e_i)=\sum\limits_{j=i}^n \frac{(-1)^{j-i}a_{1,1}^{i-2} a_{2,2}a_{1,n+1}^{j-i}}{(j-i)!} e_j, \quad  2\leq i\leq n.
 \end{equation}

 Obviously, the equality holds for $i = 2$. Let us assume that the equality holds
for $2 < i < n$, and we shall prove it for $i + 1$:
\[\begin{array}{c}
\varphi(e_{i+1})=\varphi([e_i,e_1])=[\varphi(e_i),\varphi(e_1)]=[\sum\limits_{j=i}^n \frac{(-1)^{j-i}a_{1,1}^{i-2} a_{2,2}a_{1,n+1}^{j-i}}{(j-i)!} e_j,a_{1,1}e_1]=\\
\sum\limits_{j=i}^n \frac{(-1)^{j-i}a_{1,1}^{i-1} a_{2,2}a_{1,n+1}^{j-i}}{(j-i)!} [e_j,e_1]=\sum\limits_{j=i}^{n-1} \frac{(-1)^{j-i}a_{1,1}^{i-1} a_{2,2}a_{1,n+1}^{j-i}}{(j-i)!} e_{j+1}=\\
\sum\limits_{j=i+1}^{n} \frac{(-1)^{j-i-1}a_{1,1}^{i-1} a_{2,2}a_{1,n+1}^{j-i-1}}{(j-i-1)!}e_j=\sum\limits_{j=i+1}^{n} \frac{(-1)^{j-(i+1)}a_{1,1}^{(i+1)-2} a_{2,2}a_{1,n+1}^{j-(i+1)}}{(j-(i+1))!}e_j
\end{array}\]
so, the induction proves the equalities (\ref{r1}) for any $i, 2 \leq i \leq n$.

By denoting $(a_{1,1},a_{2,2},a_{1,n+1})=(\alpha,\beta,\gamma)$ we complete the necessity condition of the proof.

Now, we shall prove the sufficiency of the theorem.

Let us $u,v\in R_1$ and denote by
\[u=\sum\limits_{i=1}^n\lambda_ie_i+\lambda_{n+1}x+\lambda_{n+2}y,
\quad v=\sum\limits_{i=1}^n\mu_ie_i+\mu_{n+1}x+\mu_{n+2}y.\]

Consider
\[\begin{array}{c}
\varphi([u,v])=\varphi\big(\big[\sum\limits_{i=1}^n\lambda_ie_i+\lambda_{n+1}x+\lambda_{n+2}y,\sum\limits_{i=1}^n\mu_ie_i+\mu_{n+1}x+\mu_{n+2}y\big]\big)=\\
\varphi\big(\lambda_1 \mu_1 e_1+ \sum\limits_{i=2}^{n-1}\mu_1\lambda_ie_{i+1}-\mu_1\lambda_{n+1}e_1+\sum\limits_{i=2}^n(i-1)\lambda_i\mu_{n+1}e_i+\sum\limits_{i=2}^n\lambda_i\mu_{n+2}e_i\big)=\\
\varphi\big((\lambda_1\mu_{n+1}-\mu_1\lambda_{n+1})e_1-\mu_1\lambda_1e_2+\sum\limits_{i=2}^n(\mu_1\lambda_{i-1}+\lambda_i((i-1)\mu_{n+1}+\mu_{n+2}))e_i\big)=\\
\alpha(\lambda_1\mu_{n+1}-\mu_1\lambda_{n+1})e_1-\sum\limits_{i=2}^n \mu_1 \lambda_1 \frac{(-1)^{i-2}\beta \gamma^{i-2}}{(i-2)!}+\\
+\sum\limits_{i=2}^n\Big((\mu_1\lambda_{i-1}+\lambda_i((i-1)\mu_{n+1}+\mu_{n+2}))\sum\limits_{j=i}^n \frac{(-1)^{j-i}\alpha^{i-2} \beta\gamma^{j-i}}{(j-i)!} e_j\Big)=\\
\alpha(\lambda_1\mu_{n+1}-\mu_1\lambda_{n+1})e_1-\sum\limits_{i=2}^n \mu_1 \lambda_1 \frac{(-1)^{i-2}\beta \gamma^{i-2}}{(i-2)!}e_i +\\
+\sum\limits_{i=2}^n \sum\limits_{j=2}^i (\mu_1\lambda_{j-1}+\lambda_j((j-1)\mu_{n+1}+\mu_{n+2})) \frac{(-1)^{i-j}\alpha^{j-2} \beta\gamma^{i-j}}{(i-j)!} e_i.
\end{array}\]

On the other hand,
\[\begin{array}{c}
[\varphi(u),\varphi(v)]=[\varphi(\sum\limits_{i=1}^n\lambda_ie_i+\lambda_{n+1}x+\lambda_{n+2}y),\varphi(\sum\limits_{i=1}^n\mu_ie_i+\mu_{n+1}x+\mu_{n+2}y)]=\\[1mm]
[(\alpha\lambda_1+\gamma\lambda_{n+1})e_1+\sum\limits_{i=2}^n\lambda_i\big(\sum\limits_{j=i}^n\frac{(-1)^{j-i}\alpha^{i-2}\beta\gamma^{j-i}}{(j-i)!}e_j\big)+\lambda_{n+1}x+\lambda_{n+2}y,\\[1mm]
(\alpha\mu_1+\gamma\mu_{n+1})e_1+\sum\limits_{i=2}^n\mu_i\big(\sum\limits_{j=i}^n\frac{(-1)^{j-i}\alpha^{i-2}\beta\gamma^{j-i}}{(j-i)!}e_j\big)+\mu_{n+1}x+\mu_{n+2}y]=\\[1mm]
[(\alpha\lambda_1+\gamma\lambda_{n+1})e_1+\sum\limits_{i=2}^n\sum\limits_{j=2}^i\lambda_j\frac{(-1)^{i-j}\alpha^{j-2}\beta\gamma^{i-j}}{(i-j)!}e_i+\lambda_{n+1}x+\lambda_{n+2}y,\\[1mm]
(\alpha\mu_1+\gamma\mu_{n+1})e_1+\sum\limits_{i=2}^n\sum\limits_{j=2}^i\mu_j \frac{(-1)^{i-j} \alpha^{j-2}\beta\gamma^{i-j}}{(i-j)!}e_i+\mu_{n+1}x+\mu_{n+2}y]=\\[1mm]
(\alpha\mu_1+\gamma\mu_{n+1})\sum\limits_{i=2}^{n-1}\sum\limits_{j=2}^i\lambda_j \frac{(-1)^{i-j} \alpha^{j-2}\beta\gamma^{i-j}}{(i-j)!}e_{i+1}-\lambda_{n+1}(\alpha \mu_1+\gamma \mu_{n+1})e_1+\\[1mm]
\mu_{n+1}(\alpha\lambda_1+\gamma\lambda_{n+1})e_1+\mu_{n+1}\sum\limits_{i=2}^n\sum\limits_{j=2}^i\lambda_j \frac{(-1)^{i-j} \alpha^{j-2}\beta \gamma^{i-j}}{(i-j)!}(i-1)e_i+\\[1mm]
\mu_{n+2}\sum\limits_{i=2}^n\sum\limits_{j=2}^i\lambda_j \frac{(-1)^{i-j} \alpha^{j-2}\beta \gamma^{i-j}}{(i-j)!}e_i=\\[1mm]
=(\lambda_1 \mu_{n+1}-\lambda_{n+1} \mu_1)\alpha e_1+\sum\limits_{i=3}^{n}\sum\limits_{j=3}^i \mu_1 \lambda_{j-1}\frac{(-1)^{i-j}\alpha^{j-2}\beta \gamma^{i-j}}{(i-j)!}e_i+\\[1mm]
\end{array}\]
\[\begin{array}{c}\sum\limits_{i=3}^{n}\sum\limits_{j=3}^i \mu_{n+1} \lambda_{j-1}\frac{(-1)^{i-j}\alpha^{j-2}\beta \gamma^{i-j+1}}{(i-j)!}e_i+
\sum\limits_{i=2}^{n}\sum\limits_{j=2}^i \mu_{n+1} \lambda_{j}\frac{(-1)^{i-j}\alpha^{j-2}\beta \gamma^{i-j}}{(i-j)!}(i-1)e_i+ \\[1mm]
+\sum\limits_{i=2}^{n}\sum\limits_{j=2}^i \mu_{n+2} \lambda_{j}\frac{(-1)^{i-j}\alpha^{j-2}\beta \gamma^{i-j}}{(i-j)!}e_i \\[1mm]
=(\lambda_1 \mu_{n+1}-\lambda_{n+1} \mu_1)\alpha e_1-\lambda_1 \mu_1 \sum\limits_{i=2}^n \frac{(-1)^{i-2} \beta \gamma^{i-2}}{(i-2)!}e_i + \\[1mm]
\sum\limits_{i=2}^{n}\sum\limits_{j=2}^i \mu_1 \lambda_{j-1}\frac{(-1)^{i-j}\alpha^{j-2}\beta \gamma^{i-j}}{(i-j)!}e_i+\\[1mm]
\sum\limits_{i=2}^{n}\sum\limits_{j=2}^i \mu_{n+1} \lambda_{j}(j-1)\frac{(-1)^{i-j}\alpha^{j-2}\beta \gamma^{i-j}}{(i-j)!}e_i+\sum\limits_{i=2}^{n}\sum\limits_{j=2}^i \mu_{n+2} \lambda_{j}\frac{(-1)^{i-j}\alpha^{j-2}\beta \gamma^{i-j}}{(i-j)!}e_i=\\[1mm]
=\alpha(\lambda_1\mu_{n+1}-\mu_1\lambda_{n+1})e_1-\sum\limits_{i=2}^n \mu_1 \lambda_1 \frac{(-1)^{i-2}\beta \gamma^{i-2}}{(i-2)!}e_i +\\
+\sum\limits_{i=2}^n \sum\limits_{j=2}^i (\mu_1\lambda_{j-1}+\lambda_j((j-1)\mu_{n+1}+\mu_{n+2})) \frac{(-1)^{i-j}\alpha^{j-2} \beta\gamma^{i-j}}{(i-j)!} e_i.
\end{array}\]

Comparing coefficients at the basis elements we obtain that
\[\varphi([u,v])=[\varphi(u),\varphi(v)]\]
and we complete the proof of theorem
\end{proof}

\begin{thm}\label{12} A linear map $\varphi:R_2\to R_2$ is a automorphism if and only if when $\varphi$ has the
following form:
\[\begin{array}{ll}
     \varphi(e_1)=\alpha e_1+\sum\limits_{i=3}^n\frac{(-1)^i\alpha\beta^{i-2}}{(i-2)!}e_i, & \\
       \varphi(e_2)=\gamma e_2, & \\
      \varphi(e_i)=\sum\limits_{j=i}^n \frac{(-1)^{j-i}\alpha^{i-1} \beta^{j-i}}{(j-i)!} e_j,&  3\leq i\leq n,\\
       \varphi(x)=\beta e_1+\sum\limits_{i=3}^n\frac{(-1)^i\beta^{i-1}}{(i-1)!}e_i+x, & \\
        \varphi(y)=\delta e_2+y, & \\
     \end{array}\]
     where \(\alpha \gamma \neq0\).
\end{thm}

\begin{proof}Let us introduce the following notations:
\[\begin{array}{ll}
    \varphi(e_i)=\sum\limits_{j=1}^na_{j,i}e_j+a_{n+1,i}x+a_{n+2,i}y, & 1\leq i\leq 2,\\
      \varphi(e_k)=\varphi([e_{k-1},e_1])=[\varphi(e_{k-1}),\varphi(e_1)],&  3\leq k\leq n,\\
      \varphi(x)=\sum\limits_{j=1}^na_{j,n+1}e_j+a_{n+1,n+1}x+a_{n+2,n+1}y,  & \\
        \varphi(y)=\sum\limits_{j=1}^na_{j,n+2}e_j+a_{n+1,n+2}x+a_{n+2,n+2}y. & \\
     \end{array}\]

Similarly, as a previous cases from the equalities
\[\varphi(e_1)=\varphi([e_1,x])=[\varphi(e_1),\varphi(x)], \quad \varphi(e_2)=\varphi([e_2,y])=[\varphi(e_2),\varphi(y)]\]\text{and} \[ \varphi(e_i)=\frac{1}{i-1}\varphi([e_i,x])=\frac{1}{i-1}[\varphi(e_i),\varphi(x)], \quad 3\leq i\leq n,\]
imply that
$\varphi(e_i)\in span(e_1,\dots,e_n)$ for $1\leq i\leq n$.

Applying  the  automorphism identity to the following pairs of elements we get further constraints
\[
\begin{array}{lllll}
\text{Identity }& & & & \text{ Constraint }\\[1mm]
\hline \hline\\[1mm]
\varphi([e_2,e_1])&=&[\varphi(e_2),\varphi(e_1)]&\Longrightarrow &a_{i,2}=0, i\notin\{2,n\}, \\[1mm]
\varphi([y,e_2])&=&[\varphi(y),\varphi(e_2)]&\Longrightarrow &a_{n,2}=0, a_{n+2,n+2}=1, \\[1mm]
\varphi([e_2,x])&=&[\varphi(e_2),\varphi(x)]&\Longrightarrow &a_{n+2,n+1}=0,\\[1mm]
\varphi([y,e_1])&=&[\varphi(y),\varphi(e_1)]&\Longrightarrow &a_{2,1}=a_{i,n+2}=0, i\notin\{2,n,n+2\}, \\[1mm]
\varphi([y,x])&=&[\varphi(y),\varphi(x)]&\Longrightarrow &a_{2,n+1}=a_{n,n+2}=0,\\[1mm]
\end{array}\]

Similar arguments for the products
\[
\begin{array}{c}
a_{1,1}e_1+\sum\limits_{i=3}^na_{i,1}e_i=\varphi(e_1)=\varphi([e_1,x])=[\varphi(e_1),\varphi(x)]=\\[1mm]
[a_{1,1}e_1+\sum\limits_{i=3}^na_{i,1}e_i,a_{1,n+1}e_1+\sum\limits_{j=3}^na_{j,n+1}e_j+a_{n+1,n+1}x]=\\[1mm]
a_{1,1}a_{n+1,n+1}e_1+(a_{1,1}a_{1,n+1}+2a_{3,1}a_{n+1,n+1})e_3+\\
+\sum\limits_{i=4}^n(a_{i-1,1}a_{1,n+1}+(i-1)a_{i,1}a_{n+1,n+1})e_i
\end{array}\]
yield $a_{n+1,n+1}=1, a_{3,1}=-a_{1,1}a_{1,n+1}$ and the next recurrent expression
\[a_{i,1}=-\frac{a_{i-1,1}a_{1,n+1}}{i-2}, \quad 4\leq i\leq n.\]

From the previous recurrent expressions we can deduce that
\[a_{i,1}=\frac{(-1)^ia_{1,1}a_{1,n+1}^{i-2}}{(i-2)!}, \quad 3\leq i\leq n.\]

Applying analogous argumentations as we used above for the product
\[0=\varphi([x,x])=[\varphi(x),\varphi(x)]\]
we obtain
\[a_{3,n+1}=-\frac{a_{1,n+1}^2}{2}, a_{i,n+1}=-\frac{a_{i-1,n+1}a_{1,n+1}}{i-1}, \quad 4\leq i\leq n,\]
which implies that
\[a_{i,n+1}=\frac{(-1)^ia_{1,n+1}^{i-1}}{(i-1)!}, \quad 3\leq i\leq n.\]

From the equality
\[\begin{array}{c}\varphi(e_3)=\varphi([e_1,e_1])=[\varphi(e_1),\varphi(e_1)]=\\[1mm]
[a_{1,1}e_1+\sum\limits_{i=3}^n\frac{(-1)^ia_{1,1}a_{1,n+2}^{i-2}}{(i-2)!}e_i,a_{1,1}e_1+
\sum\limits_{i=3}^n\frac{(-1)^ia_{1,1}a_{1,n+2}^{i-2}}{(i-2)!}e_i]=\\[1mm]
a_{1,1}^2e_3+\sum\limits_{i=3}^n\frac{(-1)^ia_{1,1}^2a_{1,n+2}^{i-2}}{(i-2)!}e_{i+1}=\sum\limits_{i=4}^n\frac{(-1)^{i-1}a_{1,1}^2a_{1,n+2}^{i-3}}{(i-3)!}e_i.
\end{array}\]
we get
\[\varphi(e_3)=\sum\limits_{i=4}^n\frac{(-1)^{i-1}a_{1,1}^2a_{1,n+2}^{i-3}}{(i-3)!}e_i.\]

With a similar induction as previous cases we can prove the following equalitiy:
\[\varphi(e_i)=\sum\limits_{j=i}^n\frac{(-1)^{j-i}a_{1,1}^{i-1}a_{1,n+1}^{j-i}}{(j-i)!}e_j, \quad 3\leq i\leq n.\]

By denoting $(a_{1,1},a_{1,n+1},a_{2,2},a_{2,n+2})=(\alpha,\beta,\gamma,\delta)$ we complete the necessity condition of the proof.

The proof of sufficiency condition of theorem is carried out by applying the similar arguments used above.

\end{proof}

\begin{thm} A linear map $\varphi:R_3\to R_3$ is a automorphism if and only if when $\varphi$ has the
following form:
\[\begin{array}{ll}
     \varphi(e_1)=\alpha e_1+\sum\limits_{i=3}^n\frac{(-1)^i\alpha\beta^{i-2}}{(i-2)!}e_i, & \\
       \varphi(e_2)=\gamma e_2, & \\
      \varphi(e_i)=\sum\limits_{j=i}^n \frac{(-1)^{j-i}\alpha^{i-1} \beta^{j-i}}{(j-i)!} e_j,&  3\leq i\leq n,\\
       \varphi(x)=\beta e_1+\sum\limits_{i=3}^n\frac{(-1)^i\beta^{i-1}}{(i-1)!}e_i+x, & \\
        \varphi(y)=y, & \\
     \end{array}\]
     where \(\alpha \gamma \neq0\).
\end{thm}

\begin{proof} The proof is similar to the proof of Theorem \ref{12}.

\end{proof}

{\small
\begin{tabular}{p{9cm}}
    Karimjnov I.A.,\\
    Andijan State University named after Z.M.Bobur, Andijan 170100, Uzbekistan. \\
    Institute of Mathematics named after V.I. Romanovsky at the Academy of Sciences of the Republic of Uzbekistan, Tashkent, 100174, Uzbekistan.\\
    e-mail: iqboli@gmail.com\\
    \\
\end{tabular}

\begin{tabular}{p{9cm}}
    Umrzaqov S.M.,\\
    Andijan State University named after Z.M.Bobur, Andijan 170100, Uzbekistan. \\
    e-mail: sardor.umrzaqov1986@gmail.com\\
    \\
\end{tabular}
}

\label{lastpage}

\end{document}